\documentclass[12pt]{amsart}%
\usepackage{amsmath}
\usepackage{graphicx}
\usepackage{amscd}
\usepackage{geometry}
\usepackage{amsfonts}
\usepackage{amssymb}%
\newtheorem{theorem}{Theorem}[section]
\theoremstyle{plain}

\newtheorem{case}{Case}

\newtheorem{conjecture}{Conjecture}
\newtheorem{corollary}[theorem]{Corollary}

\newtheorem{lemma}[theorem]{Lemma}
\newtheorem{notation}{Notation}

\newtheorem{proposition}[theorem]{Proposition}
\newtheorem{remark}{Remark}

\numberwithin{equation}{section}

\begin{document}
\title[de Groot's absolute cone conjecture ]{A solution to de Groot's absolute cone conjecture}
\author{Craig R. Guilbault }
\address{Department of Mathematical Sciences, University of Wisconsin-Milwaukee,
Milwaukee, Wisconsin 53201}
\email{craigg@uwm.edu}
\thanks{The author wishes to thank Boris Okun for a very helpful conversation.}
\date{July 1, 2005}
\subjclass{Primary 57N12, 57N13, 57N15, 57P99; Secondary 57N05}
\keywords{absolute cone, absolute suspension, homology manifold, generalized manifold,}

\begin{abstract}
A compactum $X$ is an `absolute cone' if, for each of its points $x$, the
space $X$ is homeomorphic to a cone with $x$ corresponding to the cone point.
In 1971, J. de Groot conjectured that each $n$-dimensional absolute cone is an
$n$-cell. In this paper, we give a complete solution to that conjecture. In
particular, we show that the conjecture is true for $n\leq3$ and false for
$n\geq5$. For $n=4$, the absolute cone conjecture is true if and only if the
$3$-dimensional Poincar\'{e} Conjecture is true.

\end{abstract}
\maketitle

\section{Introduction}

A compactum $X$ is an \emph{absolute suspension} if for any pair of points
$x,y\in X$, the space $X$ is homeomorphic to a suspension with $x$ and $y$
corresponding to the suspension points. Similarly, $X$ is an \emph{absolute
cone} if, for each point $x\in X$, the space $X$ is homeomorphic to a cone
with $x$ corresponding to the cone point.

At the 1971 Prague Symposium, J. de Groot \cite{Gr} made the following two conjectures:

\begin{conjecture}
\label{susp-conj}Every $n$-dimensional absolute suspension is homeomorphic to
the $n$-sphere.
\end{conjecture}

\begin{conjecture}
\label{cone-conj}Every $n$-dimensional absolute cone is homeomorphic to an $n$-cell.
\end{conjecture}

In 1974, Szyma\'{n}ski \cite{Sz} proved Conjecture \ref{susp-conj} in the
affirmative for $n=1,2$ or $3$. Later, Mitchell \cite{Mi1} reproved
Szyma\'{n}ski's results, and at the same time shed some light on higher
dimensions, by showing that every $n$-dimensional absolute suspension is an
ENR homology $n$-manifold homotopy equivalent to the $n$-sphere. Still, the
`absolute suspension conjecture' remains open for $n\geq4$.

In 2005, Nadler \cite{Na} announced a proof of Conjecture \ref{cone-conj} in
dimensions $1$ and $2$. In this paper we provide a complete solution to the
`absolute cone conjecture'. In particular, we verify Conjecture
\ref{cone-conj} for $n\leq3$ and provide counterexamples for all $n\geq5$. For
$n=4$ we show that the conjecture is equivalent to the $3$-dimensional
Poincar\'{e} Conjecture, that has recently been claimed by Perelman.

\bigskip

\section{Definitions, notation, and terminology}

\subsection{Cones}

For any topological space $L$, the \emph{cone} on $L$ is the quotient space
\[
cone\left(  L\right)  =L\times\left[  0,1\right]  /L\times\left\{  0\right\}
\text{.}%
\]
Let $q:L\times\left[  0,1\right]  \rightarrow cone(L)$ be the corresponding
quotient map. We refer to $q\left(  L\times\left\{  0\right\}  \right)  $ as
the \emph{cone point} and we view $L$ as a subspace of $cone\left(  L\right)
$ via the embedding $L\leftrightarrow L\times\left\{  1\right\}
\hookrightarrow cone\left(  L\right)  $. We refer to this copy of $L$ as the
\emph{base} of the cone. For any $\left(  z,t\right)  \in L\times\left[
0,1\right]  $, denote $q\left(  z,t\right)  $ by $t\cdot z$. Thus, $0\cdot z$
represents the cone point and $1\cdot z=z$ for all $z\in L$.

For each $z\in L$, the \emph{cone line} corresponding to $z$ is the arc
\[
I_{z}=q\left(  \left\{  z\right\}  \times\left[  0,1\right]  \right)
=\left\{  \left.  t\cdot z\ \right\vert 0\leq t\leq1\right\}  ,
\]
while the \emph{open cone line} corresponding to $z$, denoted by
$\overset{\circ}{I}_{z}$, is the set%
\[
\overset{\circ}{I}_{z}=q\left(  \left\{  z\right\}  \times(0,1\right)
)=\left\{  \left.  t\cdot z\ \right\vert 0<t<1\right\}  \text{.}%
\]

For $\varepsilon\in\left(  0,1\right)  $, the \emph{subcone of radius
}$\varepsilon$ is the set
\[
cone\left(  L,\varepsilon\right)  =q\left(  L\times\left[  0,\varepsilon
\right]  \right)  =\left\{  \left.  t\cdot z\ \right\vert z\in L\text{ and
}0\leq t\leq\varepsilon\right\}
\]
Clearly, each subcone is homeomorphic to $cone(L)$. More generally, if
$\lambda:L\rightarrow(0,1)$ is continuous, the $\lambda$\emph{-warped subcone
}is defined by
\[
cone\left(  L,\lambda\right)  =\left\{  \left.  t\cdot z\ \right\vert z\in
L\text{ and }0\leq t\leq\lambda\left(  z\right)  \right\}  \text{.}%
\]
It also is homeomorphic to $cone(L)$. In fact, the following is easy to prove.

\begin{lemma}
\label{epsilon-neighborhoods}Let $L$ be a space, $\varepsilon\in\left(
0,1\right)  $, and $\lambda:L\rightarrow(0,1)$. Then there is a homeomorphism
(in fact, an ambient isotopy) $f:cone(L)\rightarrow cone(L)$ fixed on
$L\cup\left\{  \text{cone point}\right\}  $ such that $f\left(  cone\left(
L,\varepsilon\right)  \right)  =cone\left(  L,\lambda\right)  $.
\end{lemma}

By applying the above lemma, or by a similar direct proof, we also have:

\begin{lemma}
\label{partial-homogeneity}Let $L$ be a space and suppose $t\cdot z$ and
$t^{\prime}\cdot z$ are points on the same open cone line of $cone\left(
L\right)  $. Then there is a homeomorphism (in fact, an ambient isotopy)
$f:cone(L)\rightarrow cone(L)$ fixed on $L\cup\left\{  \text{cone
point}\right\}  $ such that $f\left(  t\cdot z\right)  =t^{\prime}\cdot z$.
\end{lemma}

On occasion, we will have use for the \emph{open cone} on $L$, which we view
as a subspace of $cone\left(  L\right)  $. It is defined by%
\[
opencone\left(  L\right)  =L\times\lbrack0,1)/L\times\left\{  0\right\}  .
\]

\subsection{Suspensions and mapping cylinders}

For a topological space $L$, the \emph{suspension} of $L$ is the quotient
space
\[
susp\left(  L\right)  =L\times\left[  0,1\right]  /\left\{  L\times\left\{
0\right\}  ,L\times\left\{  1\right\}  \right\}  \text{.}%
\]
In other words, the suspension of $L$ is obtained by separately crushing out
the top and bottom levels of the product $L\times\left[  0,1\right]  $. The
images of these two sets under the quotient map are called the
\emph{suspension points }of $susp\left(  L\right)  $.

Given a map $f:L\rightarrow K$ between disjoint topological spaces, the
\emph{mapping cylinder} of $f$ is the quotient space
\[
Map\left(  f\right)  =((L\times\left[  0,1\right]  )\sqcup K)/\sim
\]
where $\sim$ is the equivalence relation on the disjoint union $(L\times
\left[  0,1\right]  )\sqcup K$ induced by the rule: $(x,0)\sim f\left(
x\right)  $ for all $x\in L$. We view $L$ and $K$ as subsets of $Map\left(
f\right)  $ via the embeddings induced by
\begin{align*}
L  &  \leftrightarrow L\times\left\{  1\right\}  \hookrightarrow
(L\times\left[  0,1\right]  )\sqcup K\text{, and}\\
K  &  \hookrightarrow(L\times\left[  0,1\right]  )\sqcup K.
\end{align*}
In addition, for each $z\in L$, the inclusion
\[
\left\{  z\right\}  \times\left[  0,1\right]  \hookrightarrow(L\times\left[
0,1\right]  )\sqcup K
\]
induces an embedding of an arc into $Map\left(  f\right)  $. We call the image
arc a \emph{cylinder line} and denote it by $E_{z}$.

Clearly, if the above range space $K$ consists of a single point, then
$Map\left(  f\right)  =cone\left(  L\right)  $. Similarly, $cone\left(
L\right)  $ can always be obtained as a quotient space of $Map\left(
f\right)  $ by crushing $K$ to a point. The following lemma will be useful
later. It allows us to view certain cones as mapping cylinders having one of
the cone lines as the range space.

\begin{lemma}
\label{cone-to-cylinder}Let $Y$ be a space and suppose $y\in Y$ has a
$k$-dimensional euclidean neighborhood $U$ in $Y$. Let $B_{0}^{k}$ and
$B_{1}^{k}$ be (tame) $k$-cell neighborhoods of $y$ lying in $U$ such that
$B_{1}^{k}\subseteq int\left(  B_{0}^{k}\right)  $. Then the pair $\left(
cone\left(  Y\right)  ,I_{y}\right)  $ is homeomorphic to $\left(  Map\left(
f\right)  ,I_{y}\right)  $ for some map $f:Y-int\left(  B_{1}^{k}\right)
\rightarrow I_{y}$. The homeomorphism may be chosen to be the identity on
$(Y-int\left(  B_{1}^{k}\right)  )\cup I_{y}$.

\begin{proof}
Choose a homeomorphism $h:B_{0}^{k}-int\left(  B_{1}^{k}\right)  \rightarrow
S^{k-1}\times\left[  0,1\right]  $ taking $\partial B_{0}^{k}$ to
$S^{k-1}\times\left\{  0\right\}  $ and $\partial B_{1}^{k}$ to $S^{k-1}%
\times\left\{  0\right\}  $. Then define $f:Y-int\left(  B_{1}^{k}\right)
\rightarrow I_{y}$ by\medskip%
\[
f\left(  x\right)  =\left\{
\begin{tabular}
[c]{ll}%
$t\cdot y$ & if $x\in B_{0}^{k}-int\left(  B_{1}^{k}\right)  $ and $h\left(
x\right)  \in S^{k-1}\times\left\{  t\right\}  $\\
$0\cdot y\quad$ & if $x\in Y-int\left(  B_{0}^{k}\right)  $%
\end{tabular}
\ \ \medskip\ \right.  .
\]
Since $f$ sends all points of $Y-int\left(  B_{0}^{k}\right)  $ to the cone
point $0\cdot y$, we may identify the `sub-mapping cylinder' $Map\left(
\left.  f\right\vert _{Y-int\left(  B_{0}^{k}\right)  }\right)  $ with the
`subcone' $cone\left(  Y-int\left(  B_{1}^{k}\right)  \right)  $. In addition,
it is easy to build a homeomorphism between the $\left(  k+1\right)  $-cell
$Map\left(  \left.  f\right\vert _{B_{0}^{k}-int\left(  B_{1}^{k}\right)
}\right)  $ and the $\left(  k+1\right)  $-cell $cone\left(  B_{0}^{k}\right)
$ taking $I_{y}$ identically onto $I_{y}$, and each cylinder line emanating
from an $x\in\partial B_{0}^{k}$ identically onto the corresponding cone line.
Fitting these pieces together yields the desired homeomorphism between
$Map\left(  f\right)  $ and $cone\left(  Y\right)  $. See Figure 1.%
\begin{figure}
[ptb]
\begin{center}
\includegraphics[
height=2.4898in,
width=5.7683in
]%
{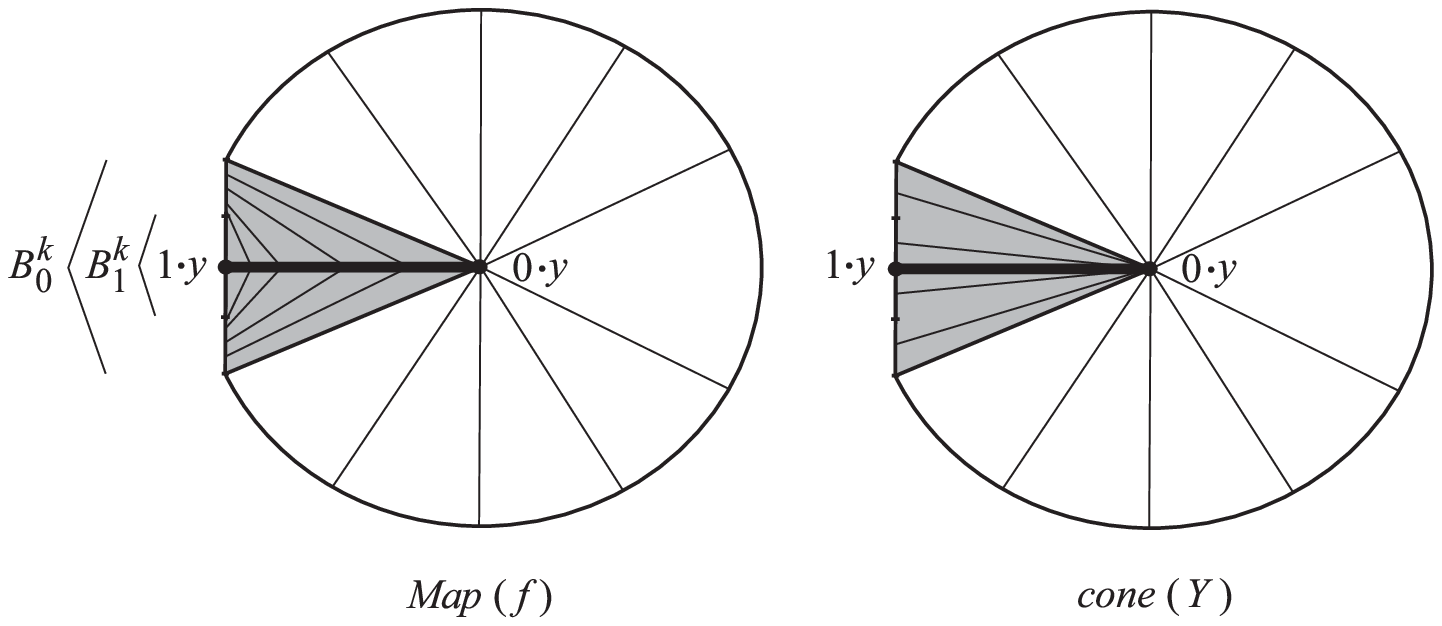}%
\caption{ }%
\end{center}
\end{figure}

\end{proof}
\end{lemma}

\subsection{ENR homology manifolds}

A space is a \emph{euclidean neighborhood retract (ENR)} if it is a retract of
some open subset of euclidean space. This is equivalent to being a
finite-dimensional separable metric ANR. A space that is a retract of
$\mathbb{R}^{n}$ (for some $n$) is called a \emph{euclidean retract (ER)}.
This is equivalent to being a contractible ENR.

A locally compact ENR $X$ is an \emph{ENR homology }$n$\emph{-manifold} if,
for every $x\in X$,%
\begin{equation}
H_{\ast}\left(  X,X-x\right)  \cong\left\{
\begin{tabular}
[c]{ll}%
$\mathbb{Z}$ & if $\ast=n$\\
$0$ & otherwise
\end{tabular}
\ \ \ \right.  \text{.} \tag{\dag$_n$}%
\end{equation}
We call $X$ an \emph{ENR homology }$n$\emph{-manifold with boundary} if, for
every $x\in X$,%
\[
H_{\ast}\left(  X,X-x\right)  \cong\left\{
\begin{tabular}
[c]{ll}%
$0$ or $\mathbb{Z}$ & if $\ast=n$\\
$0$ & otherwise
\end{tabular}
\ \ \ \right.  \text{.}%
\]
In this case, the \emph{boundary of }$X$ is the set%
\[
\partial X=\left\{  \left.  x\in X\ \right\vert H_{\ast}\left(  X,X-x\right)
\equiv0\right\}  \text{,}%
\]
and the \emph{interior of }$X$ is the set%
\[
int\left(  X\right)  =X-\partial X\text{.}%
\]
In all of the above and throughout this paper, except where stated otherwise,
homology is singular with integer coefficients.

By \cite{Mi2}, $\partial X$ is a closed subset of $X$; hence, $int\left(
X\right)  $ is an ENR homology $n$-manifold. In addition, if Borel-Moore
homology is used, $\partial X$ satisfies the algebraic condition for being a
homology $\left(  n-1\right)  $-manifold, i.e., $\partial X$ satisfies
(\dag$_{n-1}$). For ENRs, Borel-Moore homology agrees with singular homology,
so if $\partial X$ is an ENR then it is an ENR homology $\left(  n-1\right)  $-manifold.

\begin{remark}
There are interesting situations where, although $X$ is an ENR homology
manifold with boundary, $\partial X$ is not an ENR. See, for example,
\cite{AG} or \cite{Fi}. For the spaces of interest in this paper, existing
conditions will prevent this from happening.
\end{remark}

\section{Absolute cones}

Suppose a compactum $X$ is an absolute cone. For each $x\in X$, choose
$L_{x}\subseteq X$ and a homeomorphism $h_{x}:cone\left(  L_{x}\right)
\rightarrow X$ which is the identity on $L_{x}$ and sends the cone point to
$x$. We will refer to $L_{x}$ as the \emph{link} of $x$ in $X$.
(\textbf{Note.} The choice of $L_{x}$ and $h_{x}$ may not be unique; however,
for each $x$ we make a choice and stick with it.) For $\varepsilon\in\left(
0,1\right)  $ and $\lambda:L\rightarrow(0,1)$ let $N\left(  x,\varepsilon
\right)  =h_{x}\left(  cone\left(  L,\varepsilon\right)  \right)  $ and
$N\left(  x,\lambda\right)  =h_{x}\left(  cone\left(  L,\lambda\right)
\right)  $. We refer to these as the $\varepsilon$\emph{-cone neighborhood}
and the \emph{warped }$\lambda$\emph{-cone neighborhood} of $x$, respectively.
Clearly, each point of $x$ has arbitrarily small $\varepsilon$-cone neighborhoods.

In a similar vein, for any $x\in X$ and $z\in L_{x}$, let $J_{x}\left(
z\right)  $ and $\overset{\circ}{J}_{x}\left(  z\right)  $ denote
$h_{x}\left(  I_{z}\right)  $ and $h_{x}(\overset{\circ}{I}_{z})$,
respectively. We refer to these as \emph{[open] cone lines of} $X$ with
respect to $x$.

The following proposition lists several easy properties of absolute cones.

\begin{proposition}
\label{ANR}Let $X$ be a finite dimensional absolute cone, $x\in X$ and $z\in
L_{x}$. Then

\begin{enumerate}
\item $X$ is a compact ER,

\item $L_{x}$ is a compact ENR,

\item $H_{\ast}\left(  X,X-x\right)  \cong\widetilde{H}_{\ast-1}\left(
L_{x}\right)  $,

\item $L_{z}$ is contractible, and

\item $H_{\ast}\left(  X,X-z\right)  \equiv0$.
\end{enumerate}

\begin{proof}
Since each point of $X$ has arbitrarily small $\varepsilon$-cone
neighborhoods, $X$ is locally contractible; so by \cite[V.7.1]{Hu}, $X$ is an
ENR. Since $X$ is also contractible, it is an ER.

Since $L_{x}$ is a retract of its neighborhood $X-x\approx L_{x}\times
\lbrack0,1)$ in $X$, it too is an ENR \cite[III.7.7]{Hu}. Being a closed
subset of $X$, $L_{x}$ is also compact.

To prove 3), we again use $X-x\approx L_{x}\times(0,1]$. Since $X$ is
contractible, the desired isomorphisms may be obtained from the long exact
sequence for the pair $\left(  X,X-x\right)  $.

The canonical contraction of $cone\left(  L_{x}\right)  $ along cone lines
restricts to a contraction of $cone\left(  L_{x}\right)  -z$, since $z$ lies
in the base. Thus, $X-z$ is contractible. Since $X-z\approx L_{z}\times(0,1]$,
it follows that $L_{z}$ is contractible.

Assertion 5) follows from 3) and 4).
\end{proof}
\end{proposition}

The next proposition is a key ingredient in our understanding of absolute cones.

\begin{proposition}
\label{links}Let $X$ be a finite dimensional absolute cone and
\[
B_{X}=\left\{  \left.  z\in X\ \right\vert H_{\ast}\left(  X,X-z\right)
=0\right\}  .
\]
Then

\begin{enumerate}
\item $L_{x}\subseteq B_{X}$ for all $x\in X$,

\item $X-B_{X}\neq\varnothing$, and

\item For all $x\in X-B_{X}$, $L_{x}=B_{X}$.
\end{enumerate}

\begin{proof}
Assertion 1) just restates part of Proposition \ref{ANR}, while Assertion 2)
is a basic fact in dimension theory. In particular, if $\dim X=n$, then there
exists $x\in X$ such that $H_{n}\left(  X,X-x\right)  \neq0$; see, for
example, \cite[Lemma 2]{Mi2}.

To prove 3), fix $x\in X-B_{X}$ and suppose $y\in X-L_{x}$. We must show that
$y\notin B_{X}$, i.e., that $H_{\ast}\left(  X,X-y\right)  $ is non-trivial.

Choose $\varepsilon<1$ sufficiently small that $N\left(  y,\varepsilon\right)
\cap L_{x}=\varnothing$. Choose $z\in L_{y}$ such that $x$ lies on the open
cone line $\overset{\circ}{J}_{y}\left(  z\right)  $. By Lemma
\ref{partial-homogeneity}, $H_{\ast}\left(  X,X-x\right)  \cong H_{\ast
}\left(  X,X-x^{\prime}\right)  $ for all $x^{\prime}\in\overset{\circ}{J}%
_{y}\left(  z\right)  $. Therefore, $\overset{\circ}{J}_{y}\left(  z\right)
\cap L_{x}=\varnothing$. By `pushing out along' $\overset{\circ}{J}_{y}\left(
z\right)  $ we may expand the $\varepsilon$-cone neighborhood $N\left(
y,\varepsilon\right)  $ about $y$ to a warped $\lambda$-cone neighborhood
$N\left(  y,\lambda\right)  $ which contains $x$ in its interior and is
disjoint from $L_{x}$. See Figure 2.%
\begin{figure}
[ptb]
\begin{center}
\includegraphics[
height=2.7069in,
width=3.0165in
]%
{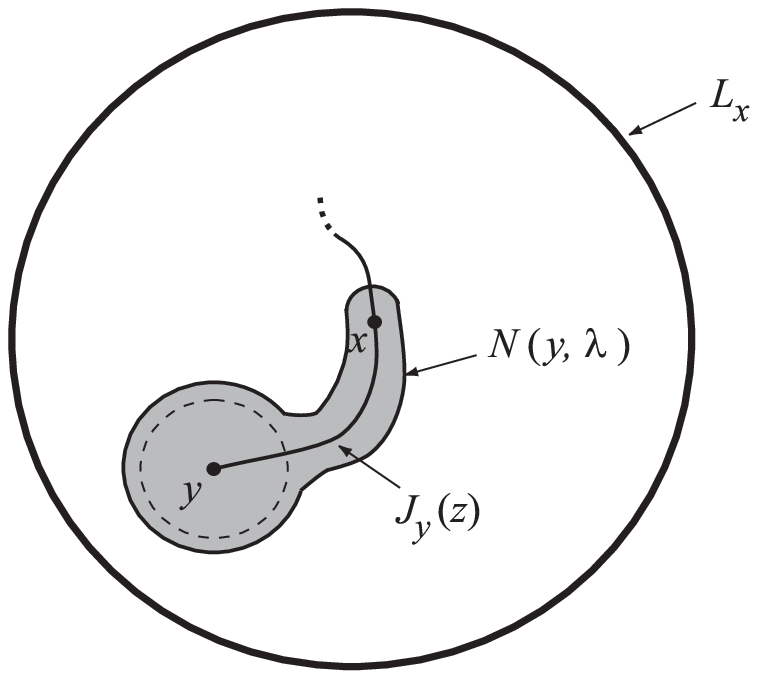}%
\caption{ }%
\end{center}
\end{figure}

Since the inclusions $\left(  X,L_{x}\right)  \hookrightarrow\left(
X,X-x\right)  $ and $\left(  X,X-N\left(  y,\lambda\right)  \right)
\hookrightarrow\left(  X,X-y\right)  $ are both homotopy equivalences of
pairs, we have inclusion induced isomorphisms%
\begin{align*}
&  H_{\ast}\left(  X,L_{x}\right)  \overset{\cong}{\longrightarrow}H_{\ast
}\left(  X,X-x\right)  \text{, and}\\
&  H_{\ast}\left(  X,X-N\left(  y,\lambda\right)  \right)  \overset{\cong
}{\longrightarrow}H_{\ast}\left(  X,X-y\right)
\end{align*}
The first of these can be factored via inclusions as such:%
\[
H_{\ast}\left(  X,L_{x}\right)  \overset{\phi}{\longrightarrow}H_{\ast}\left(
X,X-N\left(  y,\lambda\right)  \right)  \overset{\psi}{\longrightarrow}%
H_{\ast}\left(  X,X-x\right)
\]
Then $\phi$ is necessarily injective and $\psi$ surjective, so $H_{\ast
}\left(  X,X-N\left(  y,\lambda\right)  \right)  $ is non-trivial. Thus,
$H_{\ast}\left(  X,X-y\right)  \neq0$.
\end{proof}
\end{proposition}

\begin{corollary}
\label{1}For all $x\in X-B_{X}$, $H_{\ast}\left(  X,X-x\right)  \cong
\widetilde{H}_{\ast-1}\left(  B_{X}\right)  $. This homology is finitely generated.

\begin{proof}
Since $B_{X}=L_{x}$, the isomorphism follows from Proposition \ref{ANR}. Since
$B_{X}=L_{x}$ is a compact ENR, it has the homotopy type of a finite CW
complex \cite{We}; hence, finitely generated homology.
\end{proof}
\end{corollary}

\begin{corollary}
\label{2}For any $x,y\in X-B_{X}$, there exists a homeomorphism
$f:X\rightarrow X$ which is the identity on $B_{X}$ and sends $x$ to $y$.

\begin{proof}
Since $L_{x}=B_{X}=L_{y}$, we may let $f=h_{y}\circ h_{x}^{-1}$.
\end{proof}
\end{corollary}

\begin{theorem}
\label{gen-manifold-theorem}If $X$ is an $n$-dimensional absolute cone, then

\begin{enumerate}
\item $X$ is an ENR homology $n$-manifold with boundary,

\item $\partial X$ is precisely the link $L_{x}$ of any point $x\in int(X)$,

\item $\partial X$ is an ENR homology $\left(  n-1\right)  $-manifold homotopy
equivalent to $S^{n-1}$, and

\item for each $z\in\partial X$, $L_{z}$ is a contractible ENR homology
$\left(  n-1\right)  $-manifold with boundary.
\end{enumerate}

\begin{proof}
As above, let $B_{X}=\left\{  \left.  z\in X\ \right\vert H_{\ast}\left(
X,X-z\right)  \equiv0\right\}  $. By Proposition \ref{ANR} and Corollaries
\ref{1} and \ref{2}, $X-B_{X}$ is a homogeneous $n$-dimensional ENR with
finitely generated local homology. By an application of \cite{Bre} or
\cite{Bry}, $X-B_{X}$ is an ENR homology $n$-manifold. Therefore, $X$ is an
ENR homology $n$-manifold with boundary, and $\partial X=B_{X}$.

Assertion 2) is a restatement of Proposition \ref{links}.

Proposition \ref{ANR} and an application of \cite{Mi2} (as discussed earlier)
tells us that $\partial X$ is an ENR homology $\left(  n-1\right)  $-manifold.
Moreover, by Assertion 1) and Corollary \ref{1}, for any $x\in int\left(
X\right)  $,
\[
\widetilde{H}_{k-1}\left(  \partial X\right)  \cong H_{k}\left(  X,X-x\right)
\cong\left\{
\begin{tabular}
[c]{ll}%
$\mathbb{Z}$ & if $k=n$\\
$0$ & if $k\neq n$%
\end{tabular}
\ \ \ \ \ \ \ \ \ \right.
\]
Thus, $\partial X$ has the homology of $S^{n-1}$. If $n=1,2$ or $3$ it is
known that every homology $\left(  n-1\right)  $-manifold is an actual
$\left(  n-1\right)  $-manifold \cite[Ch.IX]{Wi}; and in those dimensions a
manifold is determined by its homology. Hence $\partial X$ is homeomorphic to
$S^{n-1}$. In higher dimensions we only claim a homotopy equivalence between
$S^{n-1}$ and $\partial X$. This may be obtained in the usual way if we can
show that $\partial X$ is simply connected. In particular, the Hurewicz
Theorem would then assure us that $\pi_{n-1}\left(  \partial X\right)  \cong
H_{n-1}\left(  \partial X\right)  \cong\mathbb{Z}$. A generator of $\pi
_{n-1}\left(  \partial X\right)  $ provides a degree $1$ map from $S^{n-1}$ to
$\partial X$. Since $\partial X$ is an ANR---and thus has the homotopy type of
a CW complex---a theorem of Whitehead shows that this map is a homotopy
equivalence. We hold off proving simple connectivity of $\partial X$ until
after we verify Assertion 4).

To prove 4), note that the homeomorphism $h_{z}:cone\left(  L_{z}\right)
\rightarrow X$ induces a homeomorphism of $L_{z}\times(0,1]$ onto $X-z$,
taking $L_{z}\times\left\{  1\right\}  $ onto $L_{z}$. Since $L_{z}%
\times(0,1]$ is an ENR homology $n$-manifold with boundary, $L_{z}$ is an ENR
homology $\left(  n-1\right)  $-manifold with boundary. This is an application
of \cite[Th.6]{Ra}. We have already observed (Proposition \ref{ANR}) that
$L_{z}$ is contractible.

Lastly we complete Assertion 3) by showing that $\partial X$ is simply
connected when $n\geq2$. Since the above mentioned homeomorphism $L_{z}%
\times(0,1]\rightarrow X-z$ must take (homology) boundary to boundary, it
follows that $\partial X-z$ has the structure of $L_{z}$ with an open collar
attached to $\partial L_{z}$. Thus, $\partial X-z$ is contractible; and $z$
has a neighborhood in $\partial X$ homeomorphic to a cone over $\partial
L_{z}$. Therefore, $\partial X$ may be viewed as the union of open sets
$\partial X-z$ and $U$, where $\partial X-z$ is contractible and $U$ is
homeomorphic to the open cone on $\partial L_{z}$. Since the intersection of
these sets is connected, simple connectivity follows from Van Kampen's theorem.
\end{proof}
\end{theorem}

The above proof provides some additional structure information about absolute
cones which we record as:

\begin{theorem}
If $X$ is an $n$-dimensional absolute cone, then $X$ is a contractible ENR
homology $n$-manifold with boundary and $\partial X$ is a locally conical ENR
homology $\left(  n-1\right)  $-manifold; more specifically, each point of
$\partial X$ has a neighborhood in $\partial X$ which is a cone over an ENR
homology $\left(  n-2\right)  $-manifold with the homology of $S^{n-2}$.

\begin{proof}
By the above proof, each $z\in\partial X$ has a neighborhood in $\partial X$
homeomorphic to $cone\left(  \partial L_{z}\right)  $. Since this cone lies in
$\partial X$, $\partial L_{z}$ must be an ENR. Moreover, since $L_{z}$ is an
ENR homology $\left(  n-1\right)  $-manifold with boundary, then $\partial
L_{z}$ is an ENR homology $\left(  n-2\right)  $-manifold. Lastly, since
$\partial X$ has the local homology of an $\left(  n-1\right)  $-manifold at
$z$, the homology type of $\partial L_{z}$ must be that of an $\left(
n-2\right)  $-sphere.
\end{proof}
\end{theorem}

\begin{corollary}
If $X$ is an $n$-dimensional absolute cone and $n\leq3$, then $X$ is an
$n$-cell. The same is true for $n=4$, provided the $3$-dimensional
Poincar\'{e} Conjecture is true.

\begin{proof}
If $n\leq3$, we have already observed in the proof of Theorem
\ref{gen-manifold-theorem} that $B_{X}=\partial X$ is homeomorphic to
$S^{n-1}$. Thus, $X\approx cone\left(  S^{n-1}\right)  \approx B^{n}$.

If $n=4$, we have shown that $\partial X$ is an ENR homology $3$-manifold
homotopy equivalent to $S^{3}$. In addition, we know that each point of
$\partial X$ has a neighborhood in $\partial X$ homeomorphic to the cone over
an ENR homology $2$-manifold having the homology of $S^{2}$. As above, such
ENR homology $2$-manifolds are $2$-spheres. Thus, $\partial X$ is an actual
$3$-manifold. Assuming the $3$-dimensional Poincar\'{e} Conjecture, $\partial
X\approx S^{3}$; so $X$ is a $4$-cell.
\end{proof}
\end{corollary}

\begin{remark}
At the conclusion of the next section, we will show that if there exists a
homotopy $3$-sphere $H^{3}$, not homeomorphic to $S^{3}$, then $cone\left(
H^{3}\right)  $ is a $4$-dimensional absolute cone that is not a $4$-cell.
\end{remark}

\section{Counterexamples in higher dimensions}

The main goal of this section is to construct, for all $n\geq5$,
$n$-dimensional absolute cones which are not $n$-cells. In all cases, we begin
with a non-simply connected $k$-manifold $\Sigma^{k}$ having the same
$\mathbb{Z}$-homology as $S^{k}$. Existence of such manifolds for all $k\geq3$
is well-known. Our counterexamples are obtained by first coning over
$\Sigma^{k}$, then suspending that cone. This section is primarily devoted to
proving that the resulting spaces are absolute cones, but not cells.

For completeness, we will conclude this section by showing that---if there is
a counterexample to the $3$-dimensional Poincar\'{e}---then there is also a
$4$-dimensional absolute cone that is not a $4$-cell.

We begin our construction of counterexamples in dimensions $\geq5$ with a very
general lemma.

\begin{lemma}
\label{cones-suspensions}For any space $Y$, the following are homeomorphic.

\begin{enumerate}
\item $susp\left(  cone\left(  Y\right)  \right)  $

\item $cone\left(  susp\left(  Y\right)  \right)  $

\item $cone\left(  cone\left(  Y\right)  \right)  $

\item $cone\left(  Y\right)  \times\left[  0,1\right]  $.
\end{enumerate}
\end{lemma}

Before proving this lemma, consider the map $f:cone\left(  Y\right)
\times\left[  0,1\right]  \rightarrow cone\left(  Y\right)  \times\left[
0,1\right]  $ defined by
\[
f\left(  t\cdot z,s\right)  =\left(  (st)\cdot z,s\right)  .
\]
This map takes $cone\left(  Y\right)  \times\left\{  1\right\}  $ identically
onto $cone\left(  Y\right)  \times\left\{  1\right\}  $, and each level set
$cone\left(  Y\right)  \times\left\{  s\right\}  $ to the subcone of radius
$s$ contained in $cone\left(  Y\right)  \times\left\{  s\right\}  $; finally
$cone\left(  Y\right)  \times\left\{  0\right\}  $ is taken to the cone point
of $cone\left(  Y\right)  \times\left\{  1\right\}  $. This map induces a
level-preserving embedding of $cone\left(  cone\left(  Y\right)  \right)  $
into $cone\left(  Y\right)  \times\left[  0,1\right]  $. The image of the
embedding is a particularly nice realization of $cone\left(  cone\left(
Y\right)  \right)  $ which we will denote by $\mathcal{C}^{2}\left(  Y\right)
$.

A similar map $g:cone\left(  Y\right)  \times\left[  0,1\right]  \rightarrow
cone\left(  Y\right)  \times\left[  0,1\right]  $ can be used to induce an
embedding of $susp\left(  cone\left(  Y\right)  \right)  $ into $cone\left(
Y\right)  \times\left[  0,1\right]  $. We will denote the image of that map by
$\mathcal{SC}\left(  Y\right)  $. Each of the spaces $cone\left(  Y\right)
\times\left[  0,1\right]  $, $\mathcal{C}^{2}\left(  Y\right)  $ and
$\mathcal{SC}\left(  Y\right)  $ contain as a subspace $\left\{  \text{cone
point}\right\}  \times\left[  0,1\right]  $, which we call the \emph{axis }and
denote by\emph{ }$A$. In addition, the points $($cone point$,0)$, $($cone
point,$1)$, and $($cone point$,\frac{1}{2})$ will be denoted $p_{0}$, $p_{1}$
and $p_{\ast}$, respectively. See Figure 3.%
\begin{figure}
[ptb]
\begin{center}
\includegraphics[
height=4.2704in,
width=4.2134in
]%
{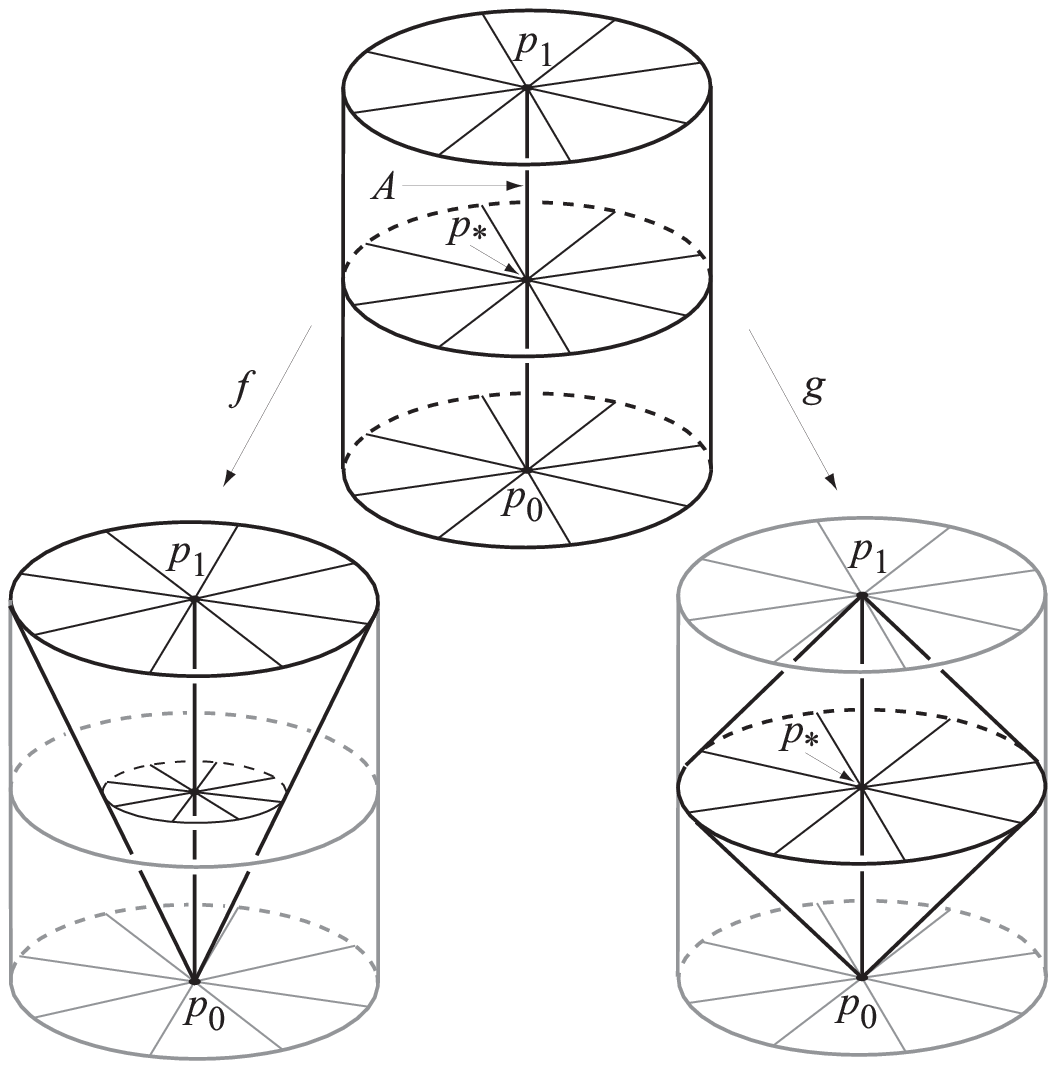}%
\caption{ }%
\end{center}
\end{figure}

\begin{proof}
[Proof of Lemma \ref{cones-suspensions}]We first show that $cone\left(
cone\left(  Y\right)  \right)  \approx cone\left(  Y\right)  \times\left[
0,1\right]  $ by illustrating a homeomorphism from $cone\left(  Y\right)
\times\left[  0,1\right]  $ to $\mathcal{C}^{2}\left(  Y\right)  $. Each cone
line $I_{y}$ of $cone\left(  Y\right)  $ determines a `square' $S_{y}%
=I_{y}\times\left[  0,1\right]  $ in $cone\left(  Y\right)  \times\left[
0,1\right]  $. Similarly, $I_{y}$ determines a `right triangle' $T_{y}$ in
$\mathcal{C}^{2}\left(  Y\right)  $ such that $T_{y}$ and $S_{y}$ have two
common sides: $I_{y}\times\left\{  1\right\}  $ and $A$. See Figure 4.%
\begin{figure}
[ptb]
\begin{center}
\includegraphics[
height=3.122in,
width=5.2641in
]%
{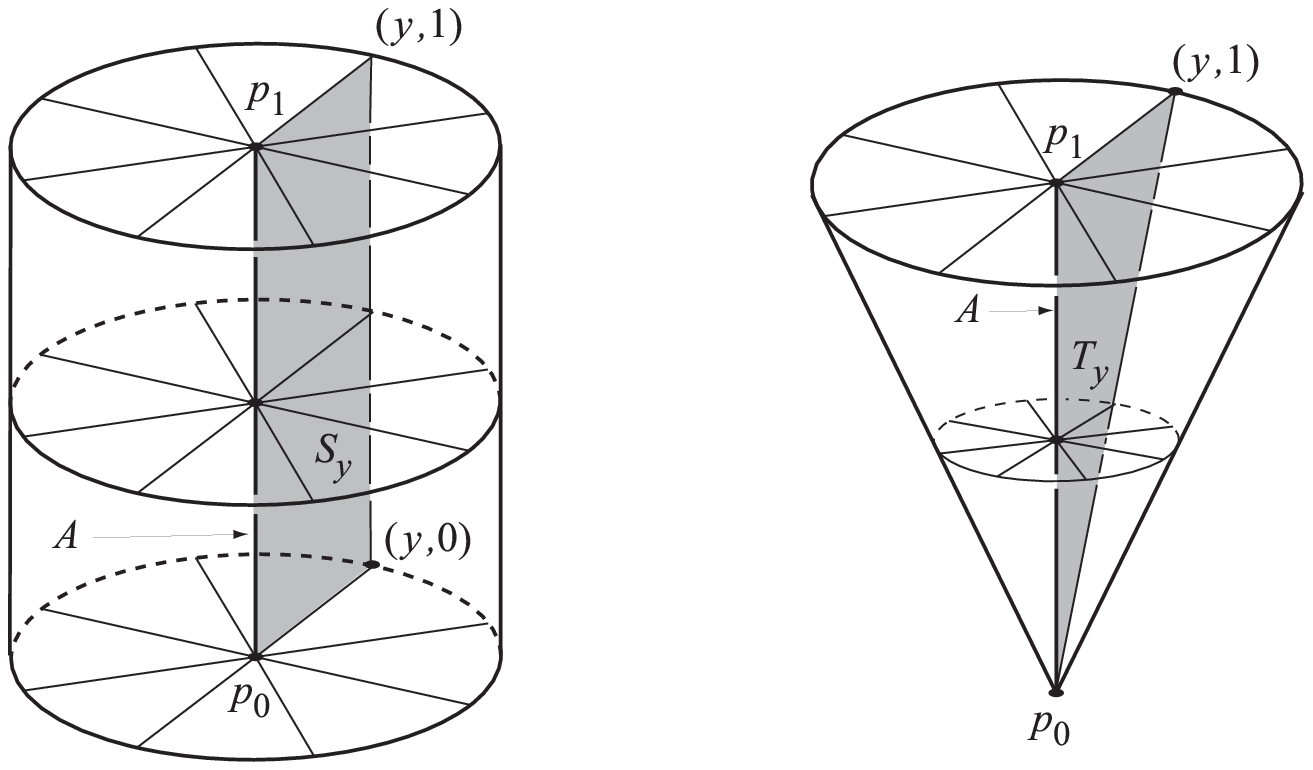}%
\caption{ }%
\end{center}
\end{figure}
For a given $y$ choose a homeomorphism
\[
k_{y}:S_{y}\rightarrow T_{y}%
\]
which is the identity on their common sides. By using the `same' homeomorphism
for each $y\in Y$, we may combine these into a single homeomorphism.%
\[
k:cone\left(  Y\right)  \times\left[  0,1\right]  \rightarrow\mathcal{C}%
^{2}\left(  Y\right)  .
\]
A similar strategy produces a homeomorphism of $cone\left(  Y\right)
\times\left[  0,1\right]  $ onto $\mathcal{SC}\left(  Y\right)  $. Thus we
have $cone\left(  cone\left(  Y\right)  \right)  \approx cone\left(  Y\right)
\times\left[  0,1\right]  \approx susp\left(  cone\left(  Y\right)  \right)  $.

Lastly, we show that $susp\left(  cone\left(  Y\right)  \right)  \approx
cone\left(  susp\left(  Y\right)  \right)  $ by observing that $\mathcal{SC}%
\left(  Y\right)  $ may be given a cone structure with base $\mathcal{S}%
\left(  Y\right)  $ and cone point $p_{\ast}$. To this end, we view
$\mathcal{SC}\left(  Y\right)  $ as the union of its `suspension triangles'
intersecting in $A$; each triangle being the suspension of a cone line. To
place a cone structure on $\mathcal{SC}\left(  Y\right)  $ we place the
obvious common cone structure on each of these triangles, with $p_{\ast}$
serving as cone point. See Figure 5.%
\begin{figure}
[ptb]
\begin{center}
\includegraphics[
height=2.674in,
width=4.9848in
]%
{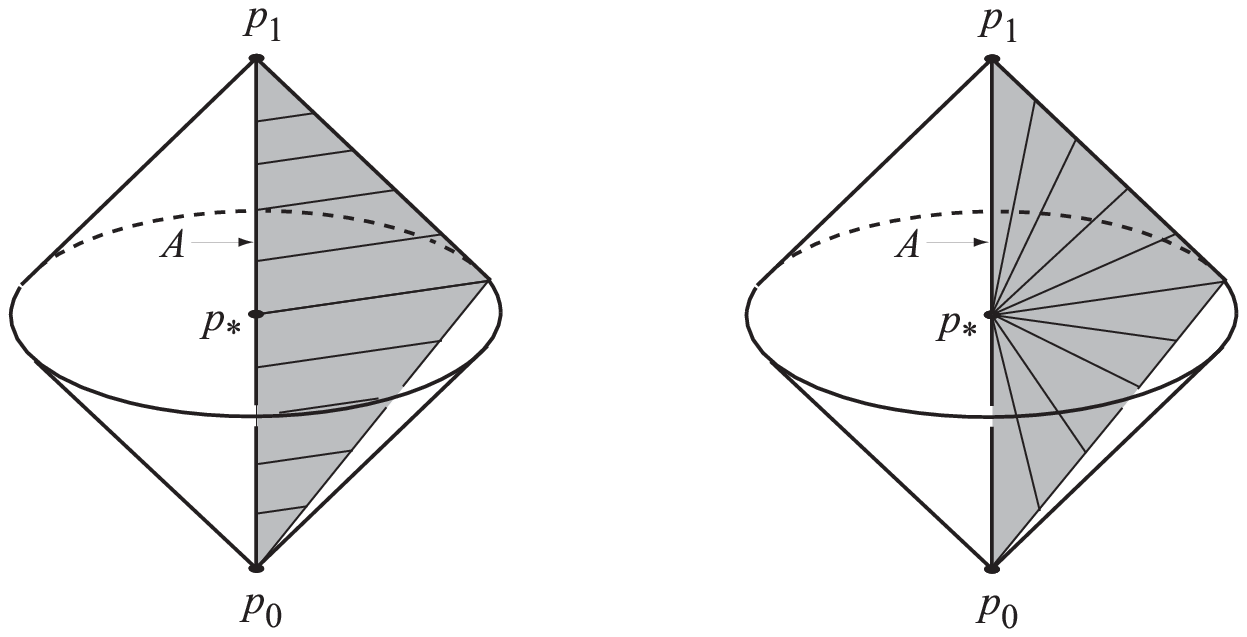}%
\caption{ }%
\end{center}
\end{figure}

\end{proof}

\begin{notation}
For convenience, we denote the each of the homeomorphic topological spaces
described by 1)-4) of Lemma \ref{cones-suspensions} by $\Theta\left(
Y\right)  $.
\end{notation}

The following lemma will be useful later. We include it now because its proof
is similar to the last part of the above argument.

\begin{lemma}
\label{cylinder-product}For any map $f:Y\rightarrow K$, the pair $\left(
Map\left(  f\right)  \times\left[  0,1\right]  ,K\times\left\{  \frac{1}%
{2}\right\}  \right)  $ is homeomorphic to $\left(  Map\left(  F\right)
,K\right)  $ for some map $F:\left(  Y\times\left[  0,1\right]  \right)
\cup(Map\left(  f\right)  \times\{0,1\})\rightarrow K$. The homeomorphism may
be chosen to take each point $\left(  k,\frac{1}{2}\right)  \in K\times
\left\{  \frac{1}{2}\right\}  $ to $k$.
\end{lemma}

\begin{proof}
To simplify matters, we assume that $f$ is surjective. If it is not, the
argument requires only minor adjustments.

By surjectivity, $Map\left(  f\right)  $ is the union of its cylinder lines
$\left\{  E_{y}\mid y\in Y\right\}  $. Hence, $Map\left(  f\right)
\times\left[  0,1\right]  $ is a union of `squares' $\left\{  S_{y}\mid y\in
Y\right\}  $, where $S_{y}$ denotes $E_{y}\times\left[  0,1\right]  $.
Furthermore, each $S_{y}$ intersects the proposed domain of our map $F$ in
three of its four boundary edges. More precisely,
\[
S_{y}\cap\left(  \left(  Y\times\left[  0,1\right]  \right)  \cup(Map\left(
f\right)  \times\{0,1\})\right)  =(\left\{  y\right\}  \times\left[
0,1\right]  )\cup\left(  E_{y}\times\left\{  0,1\right\}  \right)  .
\]
Due to its shape, we denote the right-hand side of the above equation by
$C_{y}$. Then $\left(  Y\times\left[  0,1\right]  \right)  \cup(Map\left(
f\right)  \times\{0,1\})$ is the union of the collection $\left\{  C_{y}\mid
y\in Y\right\}  $ and, for all $y,y^{\prime}\in Y$, $C_{y}\cap C_{y^{\prime}%
}=\varnothing$ unless $f(y)=f\left(  y^{\prime}\right)  $.

Define $F$ by sending each point of $C_{y}$ to $f\left(  y\right)  $. Then%
\[
F:\left(  Y\times\left[  0,1\right]  \right)  \cup(Map\left(  f\right)
\times\{0,1\})\rightarrow K
\]
is continuous, and for each $y\in Y$, the sub-mapping cylinder $Map\left(
\left.  F\right\vert _{C_{y}}\right)  $ is simply $cone(C_{y})$, with cone
point $f\left(  y\right)  \in K$. Thus, $Map\left(  F\right)  $ is a union of
the collection $\left\{  cone\left(  C_{y}\right)  \mid y\in Y\right\}  $. To
produce the desired homeomorphism from $Map\left(  f\right)  \times\left[
0,1\right]  $ to $Map\left(  F\right)  $, choose a coherent collection of
homeomorphisms from the squares $\{S_{y}\}$ making up $Map\left(  f\right)
\times\left[  0,1\right]  $ to the cones $\{cone(C_{y})\}$ making up
$Map\left(  F\right)  $. In particular, for each $y$, define a homeomorphism
$h_{y}:S_{y}\rightarrow cone(C_{y})$ which: takes $\left(  f\left(  y\right)
,\frac{1}{2}\right)  $ to the cone point $f\left(  y\right)  $, is the
identity on $C_{y}$, and is linear on segments in between these subspaces. The
union of these homeomorphisms is the desired homeomorphism from $Map\left(
f\right)  \times\left[  0,1\right]  $ to $Map\left(  F\right)  $.
\end{proof}

We are now ready to prove the main theorem of this section.

\begin{theorem}
If $\Sigma^{k}$ is a closed $k$-manifold with the same $\mathbb{Z}$-homology
as the $k$-sphere, then $\Theta\left(  \Sigma^{k}\right)  $ is a $\left(
k+2\right)  $-dimensional absolute cone. If $\Sigma^{k}$ is not simply
connected, then $\Theta\left(  \Sigma^{k}\right)  $ is not homeomorphic to an
$\left(  k+2\right)  $-cell.

\begin{proof}
We begin by observing that $\Theta\left(  \Sigma^{k}\right)  $ is ENR homology
$\left(  k+2\right)  $-manifold with boundary. To experts on homology
manifolds, this may be obvious; otherwise, proceed as follows. First note that
$cone\left(  \Sigma^{k}\right)  $ is an actual $\left(  k+1\right)  $-manifold
with boundary at all points except the cone point $p$. At that point
\begin{align*}
H_{\ast}\left(  cone\left(  \Sigma^{k}\right)  ,cone\left(  \Sigma^{k}\right)
-p\right)   &  \cong H_{\ast}\left(  cone\left(  \Sigma^{k}\right)
,\Sigma^{k}\right)  \\
&  \cong\widetilde{H}_{\ast-1}\left(  \Sigma^{k}\right)  \\
&  \cong\widetilde{H}_{\ast-1}\left(  S^{k}\right)  \cong\left\{
\begin{tabular}
[c]{ll}%
$\mathbb{Z}$ & if $\ast=k+1$\\
$0$ & otherwise
\end{tabular}
\ \ \ \ \ \right.
\end{align*}
Thus, $cone\left(  \Sigma^{k}\right)  $ is an ENR homology $\left(
k+1\right)  $-manifold with boundary; moreover, $int\left(  cone\left(
\Sigma^{k}\right)  \right)  =opencone\left(  \Sigma^{k}\right)  $. By
\cite{Ra} or straightforward calculation, $cone\left(  \Sigma^{k}\right)
\times\left[  0,1\right]  $ is an ENR homology $\left(  k+2\right)  $-manifold
with boundary, and
\[
int\left(  cone\left(  \Sigma^{k}\right)  \times\left[  0,1\right]  \right)
=opencone\left(  \Sigma^{k}\right)  \times\left(  0,1\right)  .
\]
\smallskip

\noindent\textsc{Claim 1. }$\partial\left(  \Theta\left(  \Sigma^{k}\right)
\right)  $ \emph{is not a} $\left(  k+1\right)  $\emph{-manifold. Therefore}
$\Theta\left(  \Sigma^{k}\right)  $ \emph{is not an} $\left(  n+2\right)
$\emph{-cell}.

Under the realization of $\Theta\left(  \Sigma^{k}\right)  $ as
$susp(cone\left(  \Sigma^{k}\right)  )$, the boundary is $susp\left(
\Sigma^{k}\right)  $. If $susp\left(  \Sigma^{k}\right)  $ were an actual
$\left(  k+1\right)  $-manifold, then removing a finite collection of points
would not change its fundamental group. But $susp\left(  \Sigma^{k}\right)  $
is simply connected, and $susp\left(  \Sigma^{k}\right)  -\left\{  p_{0}%
,p_{1}\right\}  \approx\Sigma^{k}\times\left(  0,1\right)  $ is not simply
connected. The claim follows.\bigskip

We now work toward showing that $\Theta\left(  \Sigma^{k}\right)  $ is an
absolute cone. For each $x\in\Theta\left(  \Sigma^{k}\right)  $, we must
identify a subspace $L_{x}$ of $\Theta\left(  \Sigma^{k}\right)  $ and a
homeomorphism of $cone\left(  L_{x}\right)  $ onto $\Theta\left(  \Sigma
^{k}\right)  $ taking the cone point to $x$. The proof splits into the
following three cases:

\begin{itemize}
\item $x\in int\left(  \Theta\left(  \Sigma^{k}\right)  \right)  ,$

\item $x=p_{0}$ or $p_{1},$

\item $x\in\partial\left(  \Theta\left(  \Sigma^{k}\right)  \right)  -\left\{
p_{0},p_{1}\right\}  $
\end{itemize}

\begin{case}
$x\in int\left(  \Theta\left(  \Sigma^{k}\right)  \right)  .$
\end{case}

In this case we know from Theorem \ref{gen-manifold-theorem} that, if $L_{x}$
exists, it must equal $\partial\left(  \Theta\left(  \Sigma^{k}\right)
\right)  $. By applying Lemma \ref{cones-suspensions} to realize
$\Theta\left(  \Sigma^{k}\right)  $ as $cone\left(  susp\left(  \Sigma
^{k}\right)  \right)  $, we see that $\Theta\left(  \Sigma^{k}\right)  $ is
indeed homeomorphic to the cone over $\partial\left(  \Theta\left(  \Sigma
^{k}\right)  \right)  $. Let $p_{\ast}\in int\left(  \Theta\left(  \Sigma
^{k}\right)  \right)  $ be the corresponding cone point and $h_{p_{\ast}%
}:cone(\partial\left(  \Theta\left(  \Sigma^{k}\right)  )\right)
\rightarrow\Theta\left(  \Sigma^{k}\right)  $ the homeomorphism. We must show
that all other elements of $int\left(  \Theta\left(  \Sigma^{k}\right)
\right)  $ can be viewed similarly. The following is the essential
ingredient.\bigskip

\noindent\textsc{Claim 2. }$int\left(  \Theta\left(  \Sigma^{k}\right)
\right)  $ \emph{is an actual} $\left(  k+2\right)  $\emph{-manifold. In fact,
}$int\left(  \Theta\left(  \Sigma^{k}\right)  \right)  \approx\mathbb{R}%
^{k+2}$.

As noted earlier, we may view $int\left(  \Theta\left(  \Sigma^{k}\right)
\right)  $ as $opencone\left(  \Sigma^{k}\right)  \times\left(  0,1\right)  $.
That this space is a $\left(  k+2\right)  $-manifold is a direct application
of work by J.W. Cannon and R.D. Edwards on the `double suspension problem'
\cite{Ca}, \cite{Ed}. A nice exposition of the relevant result may be found in
\cite[Cor. 24.3D]{Da}. Once we know that $int\left(  \Theta\left(  \Sigma
^{k}\right)  \right)  $ is a manifold, an application of \cite{St} gives us
the homeomorphism to $\mathbb{R}^{k+2}$.\bigskip

Now let $x$ be an arbitrary element of $int\left(  \Theta\left(  \Sigma
^{k}\right)  \right)  $. By a standard homogeneity argument for manifolds,
there is a homeomorphism $u:int\left(  \Theta\left(  \Sigma^{k}\right)
\right)  \rightarrow int\left(  \Theta\left(  \Sigma^{k}\right)  \right)  $
taking $p_{\ast}$ to $x$. Moreover, we may choose $u$ to be the identity
outside some compact neighborhood of $\left\{  p_{\ast},x\right\}  $. This
allows us to extend $u$ to a homeomorphism $\overline{u}$ of $\Theta\left(
\Sigma^{k}\right)  $ to itself. The homeomorphism $\overline{u}\circ
h_{p_{\ast}}:cone(\partial\left(  \Theta\left(  \Sigma^{k}\right)  )\right)
\rightarrow\Theta\left(  \Sigma^{k}\right)  $ now realizes $x$ as the cone point.

\begin{case}
$x=p_{0}$ or $p_{1}$.
\end{case}

Use Lemma \ref{cones-suspensions} to view $\Theta\left(  \Sigma^{k}\right)  $
as $cone\left(  cone\left(  \Sigma^{k}\right)  \right)  $ (or more precisely
$\mathcal{C}^{2}\left(  \Sigma^{k}\right)  $). Then $p_{0}$ corresponds to the
cone point; and the base, $cone\left(  \Sigma^{k}\right)  $, serves as
$L_{p_{0}}$. Furthermore, the realizations of $\Theta\left(  \Sigma
^{k}\right)  $ provided by 1) or 3) of Lemma \ref{cones-suspensions} reveal an
involution of $\Theta\left(  \Sigma^{k}\right)  $ interchanging $p_{0}$ and
$p_{1}$. Thus, $p_{1}$ also may be viewed as a cone point.

\begin{case}
$x\in\partial\left(  \Theta\left(  \Sigma^{k}\right)  \right)  -\left\{
p_{0},p_{1}\right\}  $.
\end{case}

By realization 1) of Lemma \ref{cones-suspensions}, $\Theta\left(  \Sigma
^{k}\right)  -\left\{  p_{1},p_{0}\right\}  \approx cone\left(  \Sigma
^{k}\right)  \times\left(  0,1\right)  $, which by another application of
Cannon-Edwards, is a $\left(  k+2\right)  $-manifold with boundary. That
boundary corresponds to $\Sigma^{k}\times\left(  0,1\right)  $. By another
standard homogeneity argument for manifolds, any two points of $\Sigma
^{k}\times\left(  0,1\right)  $ can be interchanged by a homeomorphism of
$cone\left(  \Sigma^{k}\right)  \times\left(  0,1\right)  $. This
homeomorphism can be arranged to be the identity off a compact set; and, thus,
extends to a self-homeomorphism of $susp\left(  cone\left(  \Sigma^{k}\right)
\right)  $. This means that it suffices to find a single point $x_{0}%
\in\partial\left(  \Theta\left(  \Sigma^{k}\right)  \right)  -\left\{
p_{1},p_{0}\right\}  $ at which $\Theta\left(  \Sigma^{k}\right)  $ is conical.

To this end, we return to the realization of $\Theta\left(  \Sigma^{k}\right)
$ as $cone\left(  \Sigma^{k}\right)  \times\left[  0,1\right]  $ and choose
$x_{0}\in\Sigma^{k}\times\left\{  \frac{1}{2}\right\}  $. Choose nice $k$-cell
neighborhoods $B_{1}^{k}\subseteq B_{0}^{k}$ of $x_{0}$ in $\Sigma^{k}%
\times\left\{  \frac{1}{2}\right\}  $ and, by a mild abuse of notation, let
$I_{x_{0}}$ denote the corresponding cone line as a subset of $cone(\Sigma
^{k})\times\left\{  \frac{1}{2}\right\}  $. By a combination of Lemmas
\ref{cone-to-cylinder} and \ref{cylinder-product},
\[
\left(  cone\left(  \Sigma^{k}\right)  \times\left[  0,1\right]  ,I_{x_{0}%
}\right)  \approx\left(  Map\left(  F\right)  ,I_{x_{0}}\right)
\]
for some map
\[
F:((\Sigma^{k}-int(B_{1}^{k}))\times\lbrack0,1])\cup\left(  cone\left(
\Sigma^{k}\right)  \times\left\{  0,1\right\}  \right)  \rightarrow I_{x_{0}%
}\text{.}%
\]
Crushing out the range end of the mapping cylinder yields the cone over its
domain. Hence,
\begin{equation}
cone\left(  \Sigma^{k}\right)  \times\left[  0,1\right]  /I_{x}\approx
cone\left(  ((\Sigma^{k}-int(B^{k}))\times\left[  0,1\right]  )\cup\left(
cone\left(  \Sigma^{k}\right)  \times\left\{  0,1\right\}  \right)  \right)
\text{.} \tag{\ddag}%
\end{equation}
Note that, as a subspace of the manifold with boundary $cone\left(  \Sigma
^{k}\right)  \times\left(  0,1\right)  $, $I_{x_{0}}$ is a tame arc
intersecting the boundary in a single end point. (Local tameness for this arc
is obvious at all except the interior end point of $I_{x}$. This can be seen
from the nice tubular neighborhood provided by the obvious structure of
$cone\left(  \Sigma^{k}\right)  \times\left[  0,1\right]  $. But, since our
manifold has dimension $\geq4$, $I_{x_{0}}$ cannot have a wild set consisting
of a single point (see \cite{CE} or \cite[Th.3.2.1]{Ru}). Thus $I_{x_{0}}$ is
tame in $cone\left(  \Sigma^{k}\right)  \times\left(  0,1\right)  $.)

The tameness of $I_{x_{0}}$ in $cone\left(  \Sigma^{k}\right)  \times\left[
0,1\right]  $ ensures that
\[
cone\left(  \Sigma^{k}\right)  \times\left[  0,1\right]  /I_{x_{0}}\approx
cone\left(  \Sigma^{k}\right)  \times\left[  0,1\right]  \text{.}%
\]

This homeomorphism can to chosen to take the equivalence class of $I_{x_{0}}$
to $x_{0}$. Combining this homeomorphism with (\ddag) induces the necessary
cone structure on $cone\left(  \Sigma^{k}\right)  \times\left[  0,1\right]  $
with $x_{0}$ corresponding to the cone point and
\[
L_{x_{0}}=cone\left(  ((\Sigma^{k}-int(B^{k})\times\left[  0,1\right]
)\cup\left(  cone\left(  \Sigma^{k}\right)  \times\left\{  0,1\right\}
\right)  \right)  \text{.}%
\]
\textbf{Note. }By the equivalence of tame $\left(  k+1\right)  $-cells in a
$\left(  k+1\right)  $-manifold, $L_{x_{0}}$ may be more easily visualized as
the complement of any tame open $\left(  k+1\right)  $-cell neighborhood of
$x_{0}$ in the manifold portion of $\partial\left(  \Theta\left(  \Sigma
^{k}\right)  \right)  $.
\end{proof}
\end{theorem}

Lastly, we show that the $4$-dimensional version of the absolute cone
conjecture is equivalent to the $3$-dimensional Poincar\'{e} Conjecture.

\begin{theorem}
Every $4$-dimensional absolute cone is the cone over a homotopy $3$-sphere;
moreover, if $H^{3}$ is a $3$-manifold homotopy equivalent to $S^{3}$, then
$cone\left(  H^{3}\right)  $ is an absolute cone. If $H^{3}$ is not
homeomorphic to $S^{3}$, then $cone\left(  H^{3}\right)  $ is not a $4$-cell.

\begin{proof}
We have already shown that, if $X$ is a $4$-dimensional absolute cone, then
$\partial X$ is a $3$-manifold homotopy equivalent to $S^{3}$; thus $X$ is
homeomorphic to the cone over a homotopy $3$-sphere. The last statement of the
theorem is obvious. Therefore, it remains only to show that $cone\left(
H^{3}\right)  $ is, indeed, an absolute cone. Our proof utilizes Freedman's
breakthrough work on $4$-dimensional manifolds. Aside from that application,
the proof is just a simpler version of work we have already done.\bigskip

\textsc{Claim. }$cone\left(  H^{3}\right)  $ \emph{is an absolute cone.}

By \cite[Cor.1.3]{Fr}, $H^{3}\times\mathbb{R\cong S}^{3}\times\mathbb{R}$. It
follows that $cone\left(  H^{3}\right)  $ is a $4$-manifold with
boundary---that boundary being the base $H^{3}$. By a homogeneity argument
similar to one used earlier, any point of the interior, $opencone\left(
H^{3}\right)  $, may be realized as a cone point with $H^{3}$ as its link.

If $x_{0}\in H^{3}$, use Lemma \ref{cone-to-cylinder} to realize $cone\left(
H^{3}\right)  $ as a mapping cylinder with $I_{x_{0}}$ corresponding to the
`range end'. As before, this arc is tame in the manifold $cone\left(
H^{3}\right)  $. Therefore, $cone\left(  H^{3}\right)  /I_{x_{0}}\approx
cone\left(  H^{3}\right)  $; moreover, $cone\left(  H^{3}\right)  /I_{x_{0}}$
is also homeomorphic to a cone with the equivalence class of $I_{x_{0}}$ as
its cone point (and the complement of an open $3$-cell as its base). These
homeomorphisms yield a cone structure on $cone\left(  H^{3}\right)  $ with
$x_{0}$ as the cone point.
\end{proof}
\end{theorem}

\end{document}